\documentclass[12pt, a4paper]{amsproc}

\usepackage[utf8]{inputenc}
\usepackage[russian]{babel}

\usepackage[top=20mm, bottom=20mm, outer=26mm, inner=26mm, heightrounded
]{geometry}

\usepackage[
unicode=true,
plainpages = false,
pdfpagelabels,
bookmarks=true,
bookmarksnumbered=true,
bookmarksopen=true,
breaklinks=true,
backref=false,
colorlinks=true,
linkcolor = blue,
urlcolor  = blue,
citecolor = red,
anchorcolor = green,
hyperindex = true,
hyperfigures
]{hyperref}
\hypersetup{
	pdftitle={Linear Algebra},
	pdfauthor={Vahagn H. Mikaelian},
	pdfsubject={Linear Algebra}
}

\newcommand{\Z}{{\mathbb Z}}

\newcommand{\Q}{{\mathbb Q}}
\newcommand{\Co}{{\mathbb C}}

\newtheorem{Theorem}{Теорема}[section]
\newtheorem{Corollary}[Theorem]{Следствие}
\newtheorem{Lemma}[Theorem]{Лемма}

\theoremstyle{definition}

\theoremstyle{remark}
\newtheorem{Remark}[Theorem]{Замечание}
\newtheorem{Example}[Theorem]{Пример}
\numberwithin{equation}{section}

\numberwithin{equation}{section}

\begin{document}


\makeatletter
\@namedef{subjclassname@1991}{\textbf{MSC}}
\makeatother


\title[\CYRV\cyrl\cyro\cyrzh\cyre\cyrn\cyri\cyrya, \cyrz\cyra\cyrd\cyra\cyrv\cyra\cyre\cyrm\cyrery\cyre\  \cyru\cyrn\cyri\cyrv\cyre\cyrr\cyrs\cyra\cyrl\cyrsftsn\cyrn\cyrery\cyrm\cyri\  \cyrs\cyrl\cyro\cyrv\cyra\cyrm\cyri]{\normalsize  ВЛОЖЕНИЯ, ЗАДАВАЕМЫЕ УНИВЕРСАЛЬНЫМИ СЛОВАМИ \\ 
В СВОБОДНОЙ ГРУППЕ РАНГА $2$}

\author{\footnotesize В.Г. МИКАЕЛЯН
}

\begin{abstract}
Для любой счётной группы $G = \langle\, A \mathrel{|} R\, \rangle$, заданной своими порождающими $A$ и определяющими соотношениями $R$, мы рассматриваем специфический метод вложения $G$ в некоторую $2$-порожденную группу $T$. Наше вложение явно указывает образы порождающих из $A$ в группе $T$, и явно выводит из соотношений $R$ определяющие соотношения для $T$, наследующие некоторые специальные свойства из $R$. Полученный метод может быть использован при построении явных вложений рекурсивных групп в конечно определённые группы.   
\end{abstract}

\date{\today}

\maketitle

\setcounter{tocdepth}{3}

\let\oldtocsection=\tocsection
\let\oldtocsubsection=\tocsubsection
\let\oldtocsubsubsection=\tocsubsubsection
\renewcommand{\tocsection}[2]{\hspace{-12pt}\oldtocsection{#1}{#2}}
\renewcommand{\tocsubsection}[2]{\footnotesize \hspace{6pt} \oldtocsubsection{#1}{#2}}
\renewcommand{\tocsubsubsection}[2]{ \hspace{42pt}\oldtocsubsubsection{#1}{#2}}


\section{Введение}
\noindent
Целью настоящей заметки является предложение некоторых простых правил для явного вложения любой счётной группы $G$, заданной своими порождающими и определяющими соотношениями, в такую $2$-порожденную группу $T$, что определяющие соотношения $T$ могут быть легко выведены из таковых группы $G$, и они наследуют некоторые свойства соотношений $G$, необходимые для вложений рекурсивных групп в конечно определённые группы (см. пункт~\ref{SU Preserving the structure}).

По известной теореме Хигмэна, Ноймана и Нойман любая счётная группа $G$ вложима в $2$-порожденную группу $T$ \cite{HigmanNeumannNeumann}. Этот результат, названный в учебнике Робинсона  \cite{Robinson}
\textit{``вероятно самым знаменитым из всех теорем о вложении''}, стал исходной точкой для дальнейших исследований о вложениях в $2$-порожденные группы со смежными свойствами.
Обычно такие исследования касаются случаев, когда вложение обладает дополнительным свойством (является субнормальным, вербальным и т.п.), или когда группа $T$ обладает необходимым свойством, включая свойства, наследуемые от $G$ (она разрешима, обобщенно разрешима или обобщенно нильпотентна, линейно упорядочена, аппроксимируется группами с заданным свойством, проста и т.п.). Для обзора темы см. статьи
%
\cite{Dark}--\cite{On abelian subgroups}
и процитированную в них литературу.

На самом деле, оригинальной метод вложения в \cite{HigmanNeumannNeumann} и некоторые другие конструкции вложений, процитированные выше, уже являются явными, и они уже позволяют вычислить соотношения $T$, исходя из соотношений $G$. Но нам нужен метод, который не только делает обнаружение соотношений $T$ простой, автоматизированной задачей, но и \textit{сохраняет в них некоторые закономерности}, необходимые для изучения вложений рекурсивных групп в конечно определённые группы (см. ссылки в \ref{SU Preserving the structure} ниже).

\medskip
Для задания вложения нам понадобятся следующие обозначения. 
В свободной группе
$F_2=\langle
x,y
\rangle$ ранга $2$ рассмотрим некоторые \textit{универсальные слова}:
\begin{equation}
\label{EQ definition of a_i(x,y)}
a_i(x,y) = y^{(x y^i)^{\,2}\, x^{\!-1}} 
\!\! y^{-x} 
\!\!,\quad\quad i=1,2,\ldots
\end{equation}
(здесь мы используем традиционные обозначения $x^y=y^{-1}xy$,\; $x^{-y}=(x^{-1})^y$).
Допустим, произвольная счётная группа $G$ задана своими порождающими и определяющими соотношениями: 
$$
G = \langle\, A \mathrel{|} R\, \rangle= \langle a_1, a_2,\ldots \mathrel{|} r_1, r_2,\ldots \,\rangle, 
$$
где $s$-тое соотношение 
$r_s \in R$
есть слово длины $k_s$ над буквами, например,  
$a_{i_{s,1}},\ldots,a_{i_{s,\,k_s}}$ $ \in A$.\,
Если мы заменим в $r_s$
каждый $a_{i_{s,j}}$,\; $j=1,\ldots,k_s$,\; соответствующим словом 
$a_{i_{s,j}}(x,y)$, определенным выше,  получим новое слово
\begin{equation}
\label{EQ definition of r'_s}
r'_s (x,y)=
r_s\big(a_{i_{s,1}}\!(x,y),\ldots,a_{i_{s,\,k_s}}\!(x,y)\big)
\end{equation}
над всего лишь двумя буквами $x,y$ в свободной группе $F_2$.
В этих обозначениях:

\begin{Theorem}
\label{TH universal embedding}
Для любой счётной группы $G = \langle a_1, a_2,\ldots \mathrel{|} r_1, r_2,\ldots \,\rangle $ отображение $\gamma: a_i \to a_i(x,y)$,\; $i=1,2,\ldots$\,, задает инъективное вложение $G$ в $2$-порожденную группу 
$$
T_G=\big\langle x,y 
\;\mathrel{|}\;
r'_1 (x,y),\; r'_2 (x,y),\ldots\,
\big\rangle,
$$
заданную своими соотношениями 
$r'_s (x,y)$,\;
$s=1,2,\ldots$
\end{Theorem}

Это ни что иное как другая формулировка SQ-универсальности $F_2$.
Доказательства теоремы~\ref{TH universal embedding} 
и еe модификации для групп без кручения теоремы~\ref{TH universal embedding torsion-free}
занимают пункты~\ref{SU The universal generators}--\ref{SU Some simplification for torsion free groups} ниже.
Примеры и приложения с этими вложениями можно найти в пункте~\ref{SU Examples of embeddings}. 

\medskip
Мы хотели бы подчеркнуть следующий случай, имеющий отношение к вопросу Бридсона и де ла Арпа, упомянутого в проблеме 14.10 (б) в Коуровской тетради~\cite{kourovka}:
\textit{``Найти явное вложение группы $\Q$ в конечно порождeнную группу; такая группа существует по теореме IV из  \cite{HigmanNeumannNeumann}''}.
Требуемое явное вложение  $\Q$ в $2$-порожденную группу $T$ было дано в \cite{On a Problem on Explicit Embeddings of Q} двумя методами, используя свободные конструкции и сплетения.
В настоящей заметке мы добавляем еще одно свойство: 
$\Q$ может быть явно вложено в такую группу $T$, определяющие соотношения которого конкретно перечислены. В примере~\ref{EX embedding of rational group} мы указываем явное вложение $\Q$
в $2$-порожденную группу с конкретно заданными определяющими соотношениями:
$$
T_\Q =\big\langle x,y 
\;\mathrel{|}\;
(y^s)^{(x y^s)^{\,2} x^{\!-1}} 
\!y^{-(x y^{s-1})^{\,2} x^{\!-1}}
\!\!,\;\;\;\; s=2,3\ldots
\big\rangle.
$$

Среди других недавних исследований о вложениях $\Q$ 
в конечно порожденные группы мы хотели бы кратко подчеркнуть следующие:
\cite{Darbinyan Mikaelian}
продолжает отношение \textit{линейного порядка} рациональных чисел в $\Q$ на всю $2$-порожденную группу $T$, и показывает, что вложение может быть \textit{вербальным}. В
\cite{On abelian subgroups}
замечается, что $\Q$ никогда не может быть вложено в конечно порожденную  \textit{метабелеву} группу $T$ (см. параграф 7 в \cite{On abelian subgroups} а также \cite{Finiteness conditions for soluble groups}).
А в \cite{Adian Atabekyan} строится такое явное вербальное вложение $\Q$ в $2$-порожденную группу 
$T=A_\Q(m,n)$, что  \textit{центр} $T$ совпадает с образом $\Q$, т.е. $Z(T)\cong \Q$.  
Одна из задач проблемы 14.10 (a) \cite{kourovka} -- нахождение явного вложения $\Q$ в ``естественную'' \textit{конечно определённую} группу.
\cite{The Higman operations and  embeddings} 
описывает как процедура Хигмэна может быть модифицирована для семейства групп, включающую $\Q$, таким образом, чтобы явно вложить каждую из этих групп в $2$-порождённую конечно определённую группу.
А первое прямое решение выше упомянутой задачи появилось недавно в \cite{Belk Hyde Matucci}. 
Более того, $T\!\mathcal{A}$ -- одна из конечно определённых групп, построенных в  \cite{Belk Hyde Matucci} -- является $2$-порожденной и простой.

\medskip
В заключительном пункте~\ref{SU Preserving the structure}
мы ссылаемся не основную мотивацию, приведшую нас к изучению вложений в теореме~\ref{TH universal embedding}
и в
теореме~\ref{TH universal embedding torsion-free}:
на конструктивные хигмэновские вложения  \cite{Higman Subgroups in fP groups} рекурсивных групп в конечно определённые группы.


\medskip
Когда текст этой заметки был загружен на arXiv.org, мы имели возможность обсудить тему с проф. Л.А. Бокутем, который заметил интересный параллелизм с
\cite{Shirshov 58}, где А.И. Ширшов построил элементы 
$$d_k = \Big[a \circ \big\{[\cdots (a \underbrace{\circ \, b)\circ b \cdots ]\circ b}_k\big \}\Big]\circ (a\circ b),$$
\vskip-2mm
\noindent
$k=1,2,\ldots$, 
в свободной ассоциативной алгебре $A$ с двумя порождающими $a$ и $b$. Здесь $a\circ b$ обозначает произведение  $ab-ba$ в алгебре Ли, см. детали в \S 4 в  \cite{Shirshov 58}.
Эти элементы свободно порождают свободную алгебру Ли $L(a,b)$ счётного ранга.
Они используются для определения вложения любой счётной алгебры Ли в $2$-порожденную алгебру Ли.
Множество
$\{d_k \mathrel{|} k=1,2,\ldots \}$ является ``выделяющимся'' в смысле работы \cite{Shirshov 56}.
В то же время техника нашего доказательства очень отличается от
\cite{Shirshov 58, Shirshov 56}, что делает этот параллелизм тем более любопытным.
См. также \cite{Bokut 72} и  замечание~\ref{RE reference to HNN} ниже, где мы отсылаем к оригинальной статье Хигмэна, Ноймана и Нойман
  \cite{HigmanNeumannNeumann}.

Настоящая работа поддержана совместным грантом 18RF-109 РФФИ и ГКН МОНКС РА и грантом 18T-1A306 ГКН МОНКС РА.


\section{Ссылки и некоторые вспомогательные результаты}

\noindent
За общей теоретико-групповой информации мы отсылаем к  \cite{Robinson, Kargapolov Merzljakov, Rotman}.
Если $G = \langle\, A \mathrel{|} R \,\rangle$ -- представление группы $G$ своими порождающими $A$ и определяющими соотношениями $R$, то для не пересекающегося с $A$ алфавита $B$  и для любого множества $S\subseteq F_B$ групповых слов над $B$ мы обозначим через
$\langle G, B \mathrel{|} S \,\rangle$ группу $\langle \,A \cup B \mathrel{|}  R \cup S \,\rangle$.
Если $\varphi:G\to H$ -- гомоморфизм, определенный на группе $G=\langle g_1,g_2,\ldots \rangle$  образами $\varphi(g_1)=h_1$,  $\varphi(g_2)=h_2,\ldots$\; её порождающих, мы для краткости назовем $\varphi$ гомоморфизмом,  \textit{отправляющим}
 $g_1,g_2,\ldots$ \;на\; $h_1,h_2,\ldots$

\medskip

Наши доказательства в параграфе~\ref{SE Embeddings into 2-generator groups by ``universal'' generators}
будут основаны на свободных конструкциях: 
на операции свободного произведения групп,
на свободных произведениях групп с объединённой подгруппой, и на HNN-расширениях групп с одной или с многими проходными буквами (включая и случай бесконечного количества проходных букв).

Справочную информацию об этих конструкциях можно найти в  \cite{Rotman, Lyndon Schupp, Bogopolski}.  Также мы отсылаем к нашим недавним заметкам \cite{A modified proof, Subvariety structures, The Higman operations and  embeddings} за конкретными обозначениями, используемыми и здесь для записи свободного произведения
$G*_{\varphi} H = \langle G, H \mathrel{|} a=a^{\varphi} \text{ для всех $a\in A$}\, \rangle$ групп $G$ и $H$ с подгруппами $A$ и $B$,  объединёнными с помощю изоморфизма $\varphi : A \to B$;\;
и HNN-расширения
$G*_{\varphi} t=\langle G, t \mathrel{|} a^t=a^{\varphi} \text{ для всех $a\in A$}\, \rangle$
базовой группы $G$ с помощю проходной буквы $t$, соответствующей изоморфизму
$\varphi : A \to B$ подгрупп $A,B\le G$.\;
Мы также используем HNN-расширения
$G *_{\varphi_1, \varphi_2, \ldots} (t_1, t_2, \ldots) = \langle G, t_1, t_2,\ldots \mathrel{|} a_1^{t_1}=a_1^{\varphi_1}\!\!,\; a_2^{t_2}=a_2^{\varphi_2}\!\!,\,\ldots\; \text{ для всех $a_1\in A_1$, $a_2\in A_2,\ldots$}\, \rangle$
с многими проходными буквами
$t_1, t_2, \ldots$,
соответствующими изоморфизмам
 $\varphi_1: A_1 \to B_1,\; \varphi_2: A_2 \to B_2,\ldots$ для пар подгрупп
$A_1,B_1;\; A_2,B_2; \ldots$ в $G$.

\medskip
Мы собираемся использовать специфические подгруппы в свободных произведениях с объединённой подгруппой.  
Лемма~\ref{LE subgroups in amalgamated product} 
является некоторой вариацией леммы 3.1, упомянутой на стр.~465 в \cite{Higman Subgroups in fP groups} без доказательства как \textit{``очевидная из теоремы о нормальной форме для свободных произведений с объединённой подгруппой''}.
Доказательство можно найти в пункте 2.5 в \cite{A modified proof}. 

\begin{Lemma}
\label{LE subgroups in amalgamated product}
Пусть $\Gamma = G *_\varphi \!H$ -- свободное произведение групп $G$ и $H$ с подгруппами $A \le G$ и $B \le H$, объединёнными с помощью изоморфизма 
$\varphi: A \to B$.
Если $G', H'$ -- подгруппы соответственно в $G,\, H$  такие, что для $A'=G'\cap A$ и $B'=H'\cap B$ 
 имеем 
$\varphi (A') = B'$, тогда для подгруппы $\Gamma'=\langle G',H'\rangle$ в $\Gamma$ и для ограничения $\varphi'$ гомоморфизма $\varphi$ над $A'$ мы имеем:
\begin{enumerate}
\item $\Gamma' = G'*_{\varphi'} H'$,

\item $\Gamma' \cap A = A'$ и $\Gamma' \cap B=B'$,

\item $\Gamma' \cap G = G'$ и $\Gamma' \cap H = H'$.
\end{enumerate}
\end{Lemma}

Если объединяемые подгруппы тривиальны в заданном свободном произведении с объединённой подгруппой, тогда это произведение ни что иное как обычное свободное произведение тех же групп. Применяя это наблюдение к группам $G'$ и $H'$ с тривиальными пересечениями $A'$ и $B'$, получаем:

\begin{Corollary} 
\label{CO G*H free products}
В обозначениях леммы~\ref{LE subgroups in amalgamated product}:

\begin{enumerate}
\item если $A'=G'\cap A$ и $B'=H'\cap B$ оба тривиальны,
то $\Gamma'=\langle G',H'\rangle = G'*H'$,

\item если к тому же $A'$ -- свободная группа ранга $r_1$, а $B'$ -- свободная группа ранга $r_2$, то $\Gamma'= F_r$ есть свободная группа ранга $r=r_1+r_2$.
\end{enumerate}
\end{Corollary}


\section{Вложения в $2$-порожденные группы с универсальными словами}
\label{SE Embeddings into 2-generator groups by ``universal'' generators}

\subsection{Универсальные слова в свободной группе ранга $2$}
\label{SU The universal generators}

Пусть $F=F_A$ -- свободная группа над счётным алфавитом
$A=\{a_1, a_2,\ldots\}$. Зафиксируем новую порождающую $a$, и в свободном произведении $F'=F * \langle a \rangle$ 
зададим изоморфизмы циклических подгрупп $\varphi_i:\langle a \rangle \to \langle a_i a \rangle$
отправляющие $a$ на $a_i a$,\; $i=1,2,\ldots$
Определим соответствующее HNN-расширение:
$$
P\;=\;
F' *_{\varphi_1, \varphi_2,\ldots}
(t_1, t_2,\ldots)
\;=\;
\langle
F', t_1, t_2,\ldots \mathrel{|}\;
a^{t_i}= a_i a,\;\; i=1,2,\ldots\,
\rangle.
$$
Очевидно, проходные буквы $t_i$,\, $i=1,2,\ldots$\,,\,  порождают в $P$ свободную подгруппу $X$ счётного ранга.
Легко подобрать вспомогательную $2$-порожденную группу с подгруппой, изоморфной 
$X$: в свободной группе $Y=\langle
y,z
\rangle$ элементы $t_i'=y^i z^i$\!,\, $i=1,2,\ldots$\,, свободно порождают подгруппу $X'=\langle
t_1', t_2',\ldots
\rangle$ счётного ранга. Объединяя $X$ и $X'$
с помощю изоморфизма
$\psi$, отправляющего $t_1,t_2,\ldots$ на $t'_1,t'_2,\ldots$
мы получаем группу
$$
Q\;=\;
P *_{\psi} Y \;=\;
\langle\,
P,\; y,z \mathrel{|}
t_i=t_i',\;\; i=1,2,\ldots\,
\rangle.
$$
Эта группа может быть порождена уже тремя элементами $a,y,z$, т.к. её порождающие 
$$
a_i=a_i a \cdot a^{-1} = a^{t_i}a^{-1}
= a^{t_i'}a^{-1}
\quad \text{и}\quad\;\;
t_i=t_i'=y^i z^i
\!\!,\;\;\;\;\;
i=1,2,\ldots
$$
все лежат в $\langle a,y,z \rangle$.
Согласно построению $P$, 
ни одна из нетривиальных степеней $a$ не лежит в $X$, и согласно построению $Y$,
ни одна из нетривиальных степеней
$y$ не лежит в $X'$, иными словами, пересечения 
$\langle a \rangle \cap X$ и  $\langle y \rangle\cap X'$  оба тривиальны. 
Из сказанного по следствию~\ref{CO G*H free products} вытекает, что 
$\langle
a,y \rangle$ -- свободная подгруппа ранга $1+1=2$ в $Q$. 
Используя новую проходную букву $x$ для изоморфизма $\pi:\langle  y,z\rangle \to \langle a,y\rangle$ отправляющего $y,z$ на $a,y$, мы строим:
\begin{equation}
\label{EQ formula of F_2}
F_2 =Q *_{\pi} x =
\langle
Q,x \mathrel{|}
y^x\!=a,\; z^x\!=y
\rangle
=\Big(\big((F* \langle a\rangle) *_{\varphi_1, \varphi_2,\ldots}
\!(t_1, t_2,\ldots)\big)*_{\psi} S\Big)*_{\pi} x.
\end{equation}

Это немного сбивает с толку, т.к. выше мы использовали $F_2$ для обозначения свободной группы ранга $2$ над алфавитом $\{x,y\}$. Проверим, что на самом деле  \eqref{EQ formula of F_2} совпадает с этой свободной группой. 

Во-первых, $F_2$ порожден $\{x,y\}$, т.к. $a$ и $z$ могут представлены в виде слов над $x,y$:
$$
a=y^{x}=a(x,y),\quad
z=y^{x^{-1}}\!\!=z(x,y).
$$
Используя это, мы можем представить \textit{все} выше упомянутые порождающие в виде слов над $x,y$:
\vskip-4mm
$$
t_i = t_i' =\, y^i z^i
=\; y^i x y^i x^{-1}\!
=t'_i(x,y)=t_i(x,y),
$$
\vskip-5mm
\begin{equation}
\label{EQ formula for a_(x,y)}
a_i 
= a^{t_i}a^{-1} \!
= y^{x y^i \! x  y^i \! x^{\!-1}} 
\!\! y^{-x} 
= y^{(x y^i)^{\,2}\, x^{\!-1}} 
\!\! y^{-x}
=  \;
x  (y^{-i} x^{-1})^2   
y\,
(x\, y^i)^2    x^{-2} y^{-1} \! x =a_i(x,y).
\end{equation}
Во-вторых, $F_2$ свободен над $\{x,y\}$, ибо все соотношения, требуемые нашей 
``поочередно вложенной''
свободной конструкцией
(см. правую часть в \eqref{EQ formula of F_2})
уже выполняются над словами 
$
a(x,y),\, 
a_i(x,y),\, 
t_i(x,y),\,$
$ 
t'_i(x,y),\, 
y,\, 
z(x,y),\,  
x
$
в $F_2$, как это нетрудно проверить:
\begin{equation}
\label{EQ deducng everything from x, y}
\begin{split}
& \quad \;\;  a(x,y)^{t_i(x,y)}
= (y^x)^{\, y^i x y^i x^{-1}}\!\!\!
= y^{x y^i x y^i x^{-1}} \!\!\! y^{-x}\! y^{x}=
a_i(x,y)\, a(x,y),
\\
& t'_i(x,y)=t_i(x,y)
,\quad \quad
y^x = a(x,y),\quad\quad
z(x,y)^x = 
\big(y^{x^{-1}}\big)^x = y.
\end{split}
\end{equation}
Ни одно соотношение не понадобилось для реального ``увязывания'' $x$ с $y$, так что \eqref{EQ formula of F_2} свободен над $x,y$. 

Очевидно, $\bar F=\big\langle a_1(x,y), a_2(x,y), \ldots\big\rangle$ является изоморфной копией $F$ внутри
$F_2$,
и для любого слова $r(a_{i_1},\ldots,a_{i_k})\in F$ мы имеем слово
$$r' (x,y)=
r\,\big(a_{i_1}(x,y),\ldots,a_{i_k}(x,y)\big)\in \bar F \le F_2,$$ полученное заменой каждого  $a_{i_j}$ на $a_{i_j}(x,y)$,\; $j=1,\ldots,k$.
%

\subsection{Конструкция вложения}
\label{SU The embedding construction}

Допустим, счётная группа $G$ задана как $G= \langle\, A \mathrel{|} R\, \rangle = \langle a_1, a_2,\ldots \mathrel{|} r_1, r_2,\ldots \,\rangle $ где $s$-тое соотношение
$r_s(a_{i_{s,1}},\ldots,a_{i_{s,\,k_s}}\!)$
есть слово над $k_s$ буквами, как было упомянуто во Введении.
К такому соотношению мы ставим в соответствие слово
$
r'_s (x,y)$,
определенее в $\bar F$ согласно 
\eqref{EQ definition of r'_s}. Множество подобных слов для всех $s=1,2,\ldots$ составляют подмножество   $\bar R =\big\{r'_1 (x,y), r'_2 (x,y),\ldots\big\}$ с нормальным замыканием $\bar N = \langle \bar R\rangle^{\bar F}$ в $\bar F$.
Ввиду $\bar F \cong F$ мы имеем $\bar N \cong N$ и $\bar F / \bar N \;\cong\; F/ N \;\cong\; G$.
Далее определим $\tilde N = \langle \bar R\rangle^{F_2}$ как нормальное замыкание $\bar R$ во всей группе $F_2$.

Естественный гомоморфизм $\nu: F\to G \cong F/N$ отправляет $a_i$ на $g_i=N a_i$, $i=1,2,\ldots$\,,\;
и мы можем рассмотреть 
``обновленный'' изоморфизм 
$\varphi_i:\langle a \rangle \to \langle g_i a \rangle$ циклических подгрупп в свободном произведении
$G * \langle a \rangle$
(для простоты мы не вводим новые буквы для этих $\varphi_i$).
По аналогии с \eqref{EQ formula of F_2} 
мы строим еще одну 
``поочередно вложенную'' свободную конструкцию:
\begin{equation}
\label{EQ formula of H}
H=
\Big(\big((G * \langle a \rangle) *_{\varphi_1, \varphi_2,\ldots}
\!(t_1, t_2,\ldots)\big)*_{\psi} Y\Big)*_{\pi} x,
\end{equation}
используя эти новые изоморфизмы
$\varphi_i$,
а также прежние изоморфизмы
$\psi$ и $\pi$,
примененные выше для определения $Q$ и $F_2$.

Естественный гомоморфизм $\nu$
может быть продолжен до гомоморфизма $\bar \nu$ из группы $F_2$ на $H$ потребовав, чтобы $\bar\nu$ совпадал с $\nu$ над $F$, и фиксировал каждый из остальных порождающих 
$a,t_1,t_2,\ldots,\,y,z,x$.

Оказывается, что соотношения 
$\bar R$ уже достаточны для определения группы $H$, т.к.  из \eqref{EQ formula of H} ясно, что все остальные равенства 
$a^{t_i} = g_i a$,\; $t_i=t'_i$, $y^x = a$,\; $z^x=y$ в $H$ следуют (как в \eqref{EQ deducng everything from x, y}) 
из представления порождающих через $x,y$.

Т.к. $G$ тривиально вкладывается в $H$,  подгруппа $ (\bar F   \tilde N) / \tilde N $ группы $F_2/\tilde N$ изоморфна $\bar F / \bar N\cong G$.
С другой стороны, мы имеем $ (\bar F   \tilde N) / \tilde N \cong \bar F /(\bar F \cap \tilde N )$.
Но, т.к. $\bar N$ является ядром естественного гомоморфизма из $\bar F$ в $\bar F /\bar N \cong G$, то мы имеем $\bar F \cap \tilde N \le \bar N$.
И, т.к. к тому же $\bar N \le \bar F$ и $\bar N \le \tilde F$, то имеем и $
\bar F \cap \tilde N  = \bar N$. 

Эта конструкция вместе с~\ref{SU The universal generators} доказывает следующую техническую лемму:

\begin{Lemma}
\label{LE universal elements in F_2}
Пусть $F_2=\langle x,y\rangle$ -- свободная группа ранга $2$ с элементами $a_i(x,y)$, определенными в \eqref{EQ definition of a_i(x,y)}, и порождающими подгруппу $\bar F=\big\langle a_i(x,y) \mathrel{|} i=1,2,\ldots \big\rangle$ в $F_2$.
Для любой подгруппы $N$ в $\bar F$ пусть $\bar N$ 
и $\tilde N$ -- нормальные замыкания $N$ в $\bar F$ и в  $F_2$ соответственно. 
Тогда:
$$
\bar F \cap \tilde N  = \bar N.
$$
\end{Lemma}

Теперь мы можем завершить доказательство теоремы~\ref{TH universal embedding} следующим образом. 
Пусть отображение
$\gamma$,
группа
$T_G=\big\langle x,y 
\mathrel{|}
r'_1 (x,y),\; r'_2 (x,y),\ldots\,
\big\rangle$
и соотношения
$r'_s (x,y)=
r_s\big(a_{i_{s,1}}\!(x,y),\ldots,a_{i_{s,k_s}}\!(x,y)\big)$
те же самые, что упомянуты в теореме.
Т.к. элементы $a_1(x,y), a_2(x,y), \ldots$ порождают изоморфную копию $\bar F$ группы $F$ в $F_2$, то мы имеем $G \cong \bar F / \bar N$.
Т.к. ввиду леммы~\ref{LE universal elements in F_2} $\bar N= \bar F \cap \tilde N$,  то:
$$
G\;\cong\; \bar F / \bar N \;=\; \bar F / (\bar F \cap \tilde N)
\;\cong\;
(\bar F \tilde N) /\tilde N.
$$
Но $\bar F \tilde N$ есть вся группа $F_2$ и, значит, $(\bar F \tilde N) /\tilde N$ очевидно является подгруппой в $F_2 /\tilde N \;\cong\; T$.

\medskip
Теорема~\ref{TH universal embedding} предоставляет очень легкий способ вложения счётной группы $G = \langle a_1, a_2,\ldots \mathrel{|} r_1, r_2,\ldots\, \rangle $ в $2$-порожденную группу $T_G$,
соотношения которой тривиально получаются всего лишь заменяя в 
$r_1, r_2,\ldots$ 
все буквы $a_1,a_2,\ldots$  на выражения $a_1(x,y),\; a_2(x,y),\ldots$

Заметим, что $T_G$ зависит от использованного представления $G = \langle\, A \mathrel{|} R\, \rangle$, и при ином выборе $A$ и $R$ мы могли бы построить иную $2$-порожденную группу. Тем не менее мы не хотим обозначать её через 
$T_{\langle\, A \,\mathrel{|} \,R\, \rangle}$, ибо это привело бы к громоздким обозначениям в примерах пункта
\ref{SU Examples of embeddings}.

\subsection{Некоторое упрощение для групп без кручения}
\label{SU Some simplification for torsion free groups}

Изоморфизмы $\varphi_i:\langle a \rangle \to \langle a_i a \rangle$,
отправляющие $a$ на $a_i a$, использованные \ref{SU The universal generators}, в общем случае не могут быть заменены изоморфизмами, отправляющими $a$ на $a_i$, $i=1,2,\ldots$, ибо когда $g_i \!\in\! G$ из \ref{SU The embedding construction} -- элемент  \textit{конченого} порядка, тогда $\langle a \rangle$ и $\langle g_i \rangle$  \textit{не} изоморфны, и они уже не могут быть использованы в качестве ассоциированных подгрупп.
Именно это и было причиной, по которой мы вместо них использовали $g_i a$ вместо $g_i$. Но когда
$G$ --
группа \textit{без кручения}, это препятствие устранено, и мы можем заменить
$a_i a$ на $a_i$.
Это позволяет заменить $a_i(x,y)$ в \eqref{EQ formula for a_(x,y)} более коротким словом
\begin{equation}
\label{EQ definition of bar a_i(x,y)}
\bar a_i(x,y)
= a^{t_i}
=
y^{(x y^i)^{2} x^{\!-1}}
\!\!
,\quad\quad i=1,2,\ldots
\end{equation}
Заменяя в $r_s$
каждый из $a_{i_{s,j}}$ на 
$\bar a_{i_{s,j}}(x,y)$, мы получаем другое, более короткое чем $r'_s (x,y)$ слово
$$
r''_s (x,y)=
r_s\big(\bar a_{i_{s,1}}\!(x,y),\ldots,\bar a_{i_{s,\,k_s}}\!(x,y)\big)
$$
над буквами $x,y$ в свободной группе $F_2$. Имеем следующий аналог теоремы~\ref{TH universal embedding}:

\begin{Theorem}
\label{TH universal embedding torsion-free}
Для любой счётной группы без кручения $G\, = \,\langle a_1, a_2,\ldots \mathrel{|} r_1, r_2,\ldots \,\rangle $ отображение $\gamma: a_i \to \bar a_i(x,y)$,\; $i=1,2,\ldots$\,, задает инъективное вложение $G$ в $2$-порожденную группу 
$$
T_G=\big\langle x,y 
\;\mathrel{|}\;
r''_1 (x,y),\; r''_2 (x,y),\ldots\,
\big\rangle
$$
заданную своими соотношениями
$r''_s (x,y)$,\;
$s=1,2,\ldots$
\end{Theorem}

Адаптация доказательства в  \ref{SU The universal generators} и в \ref{SU The embedding construction} для этого случая тривиальна.

\begin{Remark}
	\label{RE reference to HNN}
Читатель может сопоставить приведенные выше конструкции со страницами 252--254 в \cite{HigmanNeumannNeumann}. 
Мы использовали некоторые идеи оттуда и из \cite{Higman Subgroups in fP groups}, но наше доказательство короче, и мы произвели более короткие слова $a_i(x,y)$. 
Сравните их со словами
$
e_i= 
a^{-1} b^{-1} a\, b^{-i} a\, b^{-1} a^{-1} b^{i}a^{-1} b\, a\, b^{-i} a\, b\, a^{-1} b^i
$,
использованными в \cite{HigmanNeumannNeumann}.
К тому же у нас имеются еще более короткие слова $\bar a_i(x,y)$ 
для групп без кручения.
\end{Remark}

\subsection{Примеры явных вложений}
\label{SU Examples of embeddings}

Вот некоторые применения метода с
теоремой~\ref{TH universal embedding}
и с теоремой~\ref{TH universal embedding torsion-free}.

\begin{Example}
\label{EX embedding of free abellian into 2-generator group}
Свободная абелева группа $G=\Z^\infty$ счётного ранга может быть задана как
$$
G = \big\langle a_1, a_2,\ldots \mathrel{|} [a_k,a_l],\; k,l=1,2\ldots \big\rangle
$$ 
своими соотношениями $r_s=r_{k,l}=[a_k,a_l]$.
Т.к. $G$ -- без кручения, мы можем использовать более короткую формулу \eqref{EQ definition of bar a_i(x,y)} для отображения каждого из $a_i$ на соответствующий $\bar a_i (x,y)$. 
Это определяет вложение  $\Z^\infty$ в $2$-порожденную группу:
$$
T_{\Z^\infty} =\big\langle x,y 
\;\mathrel{|}\;
\big[
y^{(x y^k)^{\,2} x^{\!-1}} 
\!\!\!\!\!,\,\,\,
y^{(x y^l)^{\,2} x^{\!-1}} 
\big] ,\;\;\; k,l=1,2\ldots
\big\rangle.
$$
\end{Example}

\begin{Example}
\label{EX embedding of rational group}
Аддитивная группа рациональных чисел $G=\Q$ 
может быть представлена \cite{Johnson} как:
$$
G = \big\langle a_1, a_2,\ldots \mathrel{|} a_s^s=a_{s-1},\; s=2,3\ldots \big\rangle,
$$ 
где порождающее $a_i$ соответствует дробному числу $\displaystyle {1 \over i!}$ с $i=2,3\ldots$
Перепишем каждый $a_s^s=a_{s-1}$ как 
$a_s^s\,a_{s-1}^{-1}$ и используем последнее как соотношение $r_s=r_s(a_{s-1},\, a_s)$ для  $s=2,3\ldots$
Т.к. $G$ также без кручения, мы можем использовать более короткую формулу \eqref{EQ definition of bar a_i(x,y)} для отображения каждого  $a_i=\displaystyle {1 \over i!}$ на $\bar a_i (x,y)$.
После легкого упрощения
$\bar a_i (x,y)=\big(y^{(x y^i)^{\,2} x^{\!-1}} \big)^i 
\big(y^{(x y^{i-1})^{\,2} x^{\!-1}}\big)^{-1}\!\!
=\,
(y^i)^{(x y^i)^{\,2} x^{\!-1}} 
y^{-(x y^{i-1})^{\,2} x^{\!-1}}
$ мы имеем вложение $\Q$ в $2$-порожденную группу:
$$
T_\Q =\big\langle x,y 
\;\mathrel{|}\;
(y^s)^{(x y^s)^{\,2} x^{\!-1}} 
y^{-(x y^{s-1})^{\,2} x^{\!-1}}
\!\!,\;\;\;\; s=2,3\ldots
\big\rangle.
$$
\end{Example}

\begin{Example}
\label{EX embedding of Pruefer group}
Квазициклическая $p$-группа Прюфера $G=\Co_{p^\infty}$ может быть представлена как:
$$
G = \big\langle a_1, a_2,\ldots \mathrel{|}
a_1^p,\;\;\, a_{s+1}^p\!=a_s
,\;\; s=1,2\ldots \big\rangle,
$$ 
где порождающее $a_i$ 
соответствует
примитивному $p^i$-тому корню $\varepsilon_i$ из единицы \cite{Kargapolov Merzljakov}.
Т.к. это уже \textit{не} группа без кручения, мы должны использовать несколько более длинную формулу
$a_i(x,y)$ из \eqref{EQ formula for a_(x,y)} в качестве образа для $a_i$.
Для первого соотношения 
$a_1^p$ группы $G$ мы имеем новое соотношение
$a_1(x,y)^p=\big(y^{(x y)^{\,2}\, x^{\!-1}} 
\!\! y^{-x}\big)^p$.
Далее, перепишем каждый из $a_{s+1}^p=a_s$ как 
$a_{s+1}^p a_s^{-1}$\!\!,\, и используем это как соотношение $r_s=r_s(a_s,\, a_{s+1})$, где  $s=1,2\ldots$\;
Соответствующим новым соотношением будет
$$
\big(y^{(x y^{s+1})^{\,2}\, x^{\!-1}} 
\!\! y^{-x}\big)^p 
\big(y^{(x y^s)^{\,2}\, x^{\!-1}} 
\!\! y^{-x}\big)^{-1}
=
\big(y^{(x y^{s+1})^{\,2}\, x^{\!-1}} 
\!\! y^{-x}\big)^p 
y^{x} y^{-(x y^s)^{\,2}\, x^{\!-1}}\!\!.
$$
И у нас есть вложение $\Co_{p^\infty}$ в $2$-порожденную группу:
$$
T_{\Co_{p^\infty}} =\big\langle x,y 
\;\mathrel{|}\;\;
\big(y^{(x y)^{\,2}\, x^{\!-1}} 
\!\! y^{-x}\big)^p\!\!,\;\;\;\;
\big(y^{(x y^{s+1})^{\,2}\, x^{\!-1}} 
\!\! y^{-x}\big)^p 
y^{x} y^{-\,(x y^s)^{\,2}\, x^{\!-1}}
\!\!\!,\;\;\; s=1,2\ldots
\big\rangle.
$$
\end{Example}

\subsection{Использование во вложениях рекурсивных групп}
\label{SU Preserving the structure}

Основная причина, по которой нам понадобились вложения
теоремы~\ref{TH universal embedding}
и
теоремы~\ref{TH universal embedding torsion-free},
касается изучения конструктивных вложений рекурсивных групп в конечно определённые группы, т.е.  \textit{конструктивных} хигмэновских вложений \cite{Higman Subgroups in fP groups} (см. в частности ссылки на вложения $\Q$ в конечно определённые группы в связи с проблемой 14.10 (a) \cite{kourovka}, упомянутые во Введении).

Одним из шагов вложения рекурсивной группы 
$G = \langle a_1, a_2,\ldots \mathrel{|} r_1, r_2,\ldots \,\rangle$  (т.е. группы с рекурсивно перечислимыми соотношениями $r_1, r_2,\ldots$) в конечно определённую группу является предварительное вложение $G$ в $2$-порожденную группу $T=T_G$. Причем эта группа $T$ также должна иметь рекурсивно перечислимое множество соотношений. 
Простые автоматизированные вложения, которые мы построили выше, действительно сохраняют это свойство.

Более того, нам нужны вложение \textit{сохраняющие специальные особенности} соотношений. 
За подробностями отсылая к \cite{The Higman operations and  embeddings}, мы даем здесь лишь общее представление. 
Каждое соотношение $2$-порожденной группы $T$ можно закодировать с помощью некоторой последовательности целых чисел.
Это позволяет изучать рекурсивно перечислимые множества соотношений с помощью определенных наборов последовательностей целых чисел  \cite{Higman Subgroups in fP groups}. Как мы видим в 
\cite{The Higman operations and  embeddings},
вложения теоремы~\ref{TH universal embedding}
и
теоремы~\ref{TH universal embedding torsion-free} 
гарантируют некоторую тесную корреляцию между соотношениями $G$ и этими наборами последовательностей, что позволяет нам строить конструктивное вложение $G$ в конечно определенную группу. 

В частности, сравните пример~\ref{EX embedding of free abellian into 2-generator group} из этой заметки с 
примером~3.1, 
с примером~3.2 и 
с примером 4.11 со ``счетами абакус'' в 
\cite{The Higman operations and  embeddings}.

{\footnotesize
\vskip4mm
\begin{tabular}{l l}
Ф-т информатики и прикладной математики
& 
Колледж технологий и инженерии\\

Ереванский государственный университет 
&
Американский университет Армении\\

ул. Алек Манукяна 1
& 
пр. Маршала Баграмяна 40\\

Ереван 0025, Армения
&
Ереван 0019, Армения\\

Эл.адр.:
\href{mailto:v.mikaelian@gmail.com}{v.mikaelian@gmail.com},
\href{mailto:vmikaelian@ysu.am}{vmikaelian@ysu.am}
$\vphantom{b^{b^{b^{b^b}}}}$

\end{tabular}
}

\end{document}